\documentclass[12pt,Йeqno]{article}
\usepackage{cmap}
\usepackage[koi8-r]{inputenc}
\usepackage{graphicx}
\usepackage{natbib}
\usepackage{amsfonts}
\usepackage[tbtags]{amsmath}
\topskip=\baselineskip%
\newcommand{\nty}{n \to \infty}

\newcommand{\lr}{\left(}

\newcommand{\rr}{\right)}
\newcommand{\lf}{\left\{}
\newcommand{\rf}{\right\}}

\newcommand{\kn}{k_n}
\newcommand{\mn}{m_n}
\newcommand{\eel}{\end{lemma}}
\newcommand{\xkn}{X_{k_n:n}}

\newcommand{\ukn}{U_{k_n:n}}

\newcommand{\xin}{X_{i:n}}

\newcommand{\uin}{U_{i:n}}
\newcommand{\wi}{W^{(n)}_i}
\newcommand{\win}{W_{i:n}}
\newcommand{\xian}{\xi_{\alpha_n}}
\newcommand{\qa}{\xi_{\alpha_n}}
\newcommand{\qb}{\xi_{1-\beta_n}}

\newcommand{\al}{\alpha}
\newcommand{\be}{\beta}
\newcommand{\si}{\sigma}
\newcommand{\an}{\alpha_n}
\newcommand{\bn}{\beta_n}
\newcommand{\xibn}{\xi_{1-\beta_n}}

\newcommand{\mau}{M_{\alpha_n}}
\newcommand{\eps}{\varepsilon}

\newcommand{\siw}{\sigma_{W,n}}
\newcommand{\tn}{T_{n}}
\newcommand{\wn}{W_{n}}
\newcommand{\dn}{\delta_n}

\numberwithin{equation}{section}

\newtheorem{theorem}{Theorem}[section]

\newtheorem{lemma}{Lemma}[section]

\newtheorem{remark}{Remark}[section]

\title{{\Huge \sf Cram\'{e}r type moderate deviations for intermediate trimmed means}
\author{\large Nadezhda V. Gribkova}
\date{
{\small \it
\centerline{St.Petersburg~State~University,~Mathematics~and~Mechanics~Faculty,}
\centerline{199034,~ Universitetskaya nab.~7/9, St.~Petersburg, Russia}
}}}
\begin{document}
\maketitle 
\begin{abstract}
{
In this article we establish
Cram\'{e}r type moderate deviation results for (intermediate) trimmed means
$T_n=n^{-1} \sum_{i=\kn+1}^{n-\mn}\xin$, where $\xin$ --  the order statistics corresponding to the first $n$  observations of a~sequence  $X_1,X_2,\dots $ of i.i.d random variables with $df$ $F$.
We consider two cases of intermediate and heavy  trimming.   In the former case, when $\max(\an,\bn)\to 0$ ($\an=\kn/n$, $\bn=\mn/n$) and  $\min(\kn,\mn)\to\infty$ as $\nty$,  we  
obtain our results under a~natural moment condition and a~mild condition on the rate at which  $\an$ and  $\bn$ tend to zero. In the latter case we do not impose any  moment conditions on $F$; instead, we require some smoothness  of  $F^{-1}$ in an~open set containing the  limit points of the trimming  sequences $\an$,  $1-\bn$.
}
\end{abstract}

\\[3mm] \noindent{\bf Mathematics Subject Classification (2010):}   62G30, 60F10,  60F05,   62E20,  62G35.\\[3mm]


\section{Introduction and main results}
\label{imtro}
Let  $X_1,X_2,\dots $ be a~sequence of independent identically
distributed (i.i.d.) real-valued 
random variables
(r.v.'s) with common distribution function ($df$) $F$, and for each
integer $n \geq 1$ let \ $X_{1:n}\le \dots \le X_{n:n}$ denote the
order statistics based on the sample $X_1,\dots ,X_n$. Introduce
the left-continuous inverse function $F^{-1}$ defined as
$F^{-1}(u)= \inf \{ x: F(x) \ge u \}$, \ $0<u\le 1$, \
$F^{-1}(0)=F^{-1}(0^+)$, and let $F_n$ and $F_n^{-1}$ denote the
empirical $df$ and its inverse respectively.

Consider the trimmed mean given by
\begin{equation}
\label{tn}
 T_n=\frac 1n \sum_{i=\kn+1}^{n-\mn}\xin
=\int_{\an}^{1-\bn}F_n^{-1}(u)\, d\, u,
\end{equation}
where $k_n, m_n$ are two sequences of integers such that $0\leq k_n<n-m_n\leq n$,
$\an=k_n/n$, \ $\bn=m_n/n$. It will assumed throughout this paper that
\begin{equation}
\label{kn_inf}
\kn\wedge \mn \to\infty
\end{equation}
as $\nty$. Here and in the sequel,  $a\wedge b:=\min(a,b,)$, $a\vee b:=\max(a,b)$. 

The  asymptotic properties of trimmed and intermediate  (when $\an \vee \bn \to 0$ along with~\eqref{kn_inf}) trimmed sums and means were investigated by many authors.
In particular in \cite{cso_hae_mas:1988}
a~necessary and sufficient condition for the existence of $\{a_n\}$, $\{b_n\}$ such that
the distribution of the properly normalized trimmed mean  $a_n^{-1}(T_n-b_n)$ tends to the standard normal law was obtained, and (using a~different approach than in \cite{cso_hae_mas:1988}) \cite{grif_pru:89} derived an~equivalent $iff$ condition for asymptotic normality of $T_n$. Both in \cite{cso_hae_mas:1988} and in ~\cite{grif_pru:89} a~classical result by \cite{stigler:1973} for
the trimmed mean with fixed trimming percentages was extended to the case that
the fraction of trimming data is vanishing when $n$ gets large.
The second order asymptotic properties (Berry-Esseen  type bounds and one-term Edgeworth  expansions) for  (intermediate)  trimmed means were established in \cite{gri_helm:2006,gri_helm:2007, gri_helm:2014}; \cite{gri:2013}. Various aspects of the bootstrap for this kind of statistics were studied, e.g., in \cite{hall_padm:1992, deh_mas_shor:1993, gri_helm:2007,gri_helm:2011} (see also the references therein).

The trimmed sums represent a~subclass of $L$-statistics. A~number of highly sharp results on Cram\'{e}r type large and moderate  deviations for $L$-statistics with smooth on $(0,1)$ weight functions -- that are not applicable  for the trimmed sum because of the discontinuity of its weights -- was obtained by \cite{vand_verav:1982, ben_zit:1990, alesk:1991}. For the case of heavy truncated $L$-statistics, when the weight function is zero outside some interval $[\al,\be]\subset (0,1)$), a~result on Cram\'{e}r type large deviations was  first established by \cite{call_vand_ver:1982};  more recently, the latter result was  strengthened in  \cite{gri:2016}, where a~different approach than in \cite{call_vand_ver:1982} was  proposed and implemented.
  

The aim of this article is to investigate  Cram\'{e}r type moderate deviations 
for intermediate
trimmed means.  To the best of our knowledge, this subject  has not been studied at all. The case of heavy trimming will be also considered.

Define  the population trimmed mean and variance functions
\begin{equation}
\label{1.1} \begin{split}
\mu (u,1-v)&=\int_u^{1-v} F^{-1}(s) \,d\, s  \ \ ,\\
\si^2 (u,1-v)&=\int_u^{1-v}\int_u^{1-v} (s\wedge t - st)\, d\,
F^{-1}(s)\, d\, F^{-1}(t),
\end{split}
\end{equation}
where $0\le u<1-v\le 1$. Note that $\si^2 (0,1)={Var}(X_1)$ whenever $\textbf{E}X_1^2$ is finite. Here and in the sequel, we use the convention that $\int_a^b=\int_{[a,b)}$ when integrating with respect to the left
continuous integrator $F^{-1}$.

Define the $\nu$-th quantile of $F$ by $\xi_{\nu}=F^{-1}(\nu)$,
$0<\nu<1$, and let $\wi$, \ $i=1,\dots,n$, denote the $X_i$
Winsorized outside of $(\qa,\qb]$. In other words
\begin{equation}
\label{wi}
\wi 
=\left\{
\begin{array}{ll}
\xian,& X_i\le \xian, \\
X_i,& \xian < X_i \le \xibn,\\
\xibn,& \xibn < X_i .
\end{array}
\right.
\end{equation}

To normalize $\tn$, we define two sequences
\begin{equation}
\label{mu_si}
\mu_n=\mu (\an,1-\bn),\qquad \siw^2={Var}(\wi).
\end{equation}
Note that $\siw^2=\si^2 (\an,1-\bn)$ and that $\mu_n$ and $\siw$ are suitable location and scale parameters for $\tn$ when
establishing its asymptotic normality (cf. Cs\"org\H{o} et al.~(1988)). We will suppose throughout
this article that $\liminf_{\nty}\siw>0$ (i.e. that $\xian
\neq \xibn$ for all sufficiently large $n$).

Let $\Phi$ denote the standard normal distribution function.
 Here is our first result on  Cram\'{e}r  type moderate deviations for the intermediate trimmed mean.
\begin{theorem}
\label{thm1}
Suppose that {\em $\textbf{E}|X_1|^p<\infty$ } for some $p>c^2+2$ $(c>0)$. In addition, assume that
\begin{equation}
\label{c_kn}
\frac{\log n}{\kn\wedge\mn}\to 0 \qquad \qquad \ \ \
\end{equation}
as $\nty$, and that
\begin{equation}
\label{c_an}
\an\vee\bn=O\bigl((\log n)^{-\gamma}\bigr),
\end{equation}
 for some $\gamma>2p/(p-2)$, as $\nty$. \
Then {\em
\begin{equation}
\label{md_1}
\begin{split}
\textbf{P}\Bigl(\frac{\sqrt{n}(\tn-\mu_n)}{\siw}>x\Bigr) &= [1-\Phi(x)](1+o(1)),\\
\textbf{P}\Bigl(\frac{\sqrt{n}(\tn-\mu_n)}{\siw} < -x\Bigr)&=\Phi(-x)(1+o(1)),
\end{split}
\end{equation}
}
as \,$\nty$, uniformly in the range $-A\leq x\leq c\sqrt{\log n}$ \  $(A>0)$.
\end{theorem}

We relegate the proof of Theorem~\ref{thm1} and other main results  to Section~\ref{proof}.

Let us discuss the conditions of Theorem~\ref{thm1}. Since our results concern  the relative error in {CLT} for $\tn$, certainly under our assumptions the $iff$  condition of the asymptotic normality of $\tn$
given in \cite{cso_hae_mas:1988} is  satisfied (see the proofs in Section~\ref{proof}). This condition is as follows: for every $t\in \mathbb{R}$
\begin{equation}
\label{q_0} \Psi_{j,n}(t)\to 0,  \ \ j=1,2,
\end{equation}
as $\nty$, where $\Psi_{j,n}$, $j=1,2$, are the auxiliary functions defined in \cite[page 674]{cso_hae_mas:1988}
which correspond to the trimming of the $k_n$
smallest and $\mn$ largest observations respectively. Consider the first of these functions (the second one is defined similarly). In our notation it is equal to
\begin{equation}
\label{psi}
\Psi_{1,n}(t)= \begin{cases} \frac{\an^{1/2}}{\si_{W_{(n)}}}\lf
F^{-1}\Bigl( \an +t\sqrt{\frac{\an}{n}}\Bigr) -F^{-1}\Bigl(\an\Bigr)\rf\ , \\
\qquad \qquad \qquad \qquad \qquad \qquad \ \qquad |c|\le \frac 12 \sqrt{\an n}\ , \\
\Psi_{1,n}\bigl(- \frac 12 \sqrt{\an n}\bigr), \ \ \ \ \qquad -\infty < t <-\frac 12 \sqrt{\an n}\ ,\\
\Psi_{1,n}\bigl(\frac 12 \sqrt{\an n}\bigr), \ \ \ \ \  \frac 12
\sqrt{\an n}<t<\infty \ .
\end{cases}
\end{equation}

If  $\an\vee\bn\to 0$ and condition~\eqref{c_an} holds, then our moment assumption, i.e., $\textbf{E}|X_1|^p<\infty$ for some $p>c^2+2$ $(c>0)$, implies~\eqref{q_0} (see the proof of Theorem~\ref{thm1}). In contrast, if $\an\vee\bn$ does not tend to zero, then (as it follows from~\eqref{psi})
the convergence in~\eqref{q_0} can happen only under some additional smoothness condition on $F^{-1}$ (cf.~condition~\eqref{cgh} in Theorem~\ref{thm3}). Thus, assumption~\eqref{c_an} in Theorem~\ref{thm1} means  that  $\an\vee\bn$ converges to zero fast enough, otherwise the trimming would be occur close to the central region,
where a~smoothness condition on $F^{-1}$ is  required, even  for  the asymptotic normality only.

As is known,  the intermediate trimmed mean $\tn$  can serve as a~consistent and robust estimator for  $\textbf{E}X_1$ (whenever it exists), and the results on large and moderate deviations for $\tn$ can be helpful  to construct  confidence intervals for the expectation of $X_1$. In particular  the question of whether is possible to replace $\mu_n$ in~\eqref{md_1} by $\textbf{E}X_1$ is of some practical interest. Our next result concerns the properties of the first two moments of $\tn$ and the possibility of replacing the normalizing sequences in~\eqref{md_1}.
\begin{theorem}
\label{thm2}
Suppose that the conditions of Theorem{\em~\ref{thm1}} hold true. Then {\em
\begin{equation}
\label{t2_1}
n^{1/2} (\textbf{E}\tn-\mu_n)=o\Bigl((\log n)^{-1}\Bigr),\qquad  \qquad
\end{equation}
\begin{equation}
\label{t2_2}
\qquad \qquad  \frac{\siw}{\si}= 1+ o\bigl( (\log n)^{-2}\bigr), \qquad
\end{equation}
\begin{equation}
\label{t2_3}
\quad\frac{\sqrt{{Var}(\tn)}}{\siw/\sqrt{n}}= 1+ o\bigl((\log n)^{-1}\bigr), \qquad
\end{equation}
}
as \,$\nty$. Moreover, $\mu_n$ and $\siw$ in relations~\eqref{md_1} can be replaced respectively by {\em $\textbf{E}\tn$ } and
$\si$ or {\em $\sqrt{n {Var}(\tn)}$}, without affecting the result.

Furthermore, if in addition  
\begin{equation}
\label{c_an2}
\an\vee\bn=o\bigl[(n\log n)^{-\frac{p}{2(p-1)}}\bigr],
\end{equation}
then  {\em
\begin{equation}
\label{t2_4}
n^{1/2} (\textbf{E}X_1-\mu_n)=o\Bigl((\log n)^{-1/2}\Bigr),
\end{equation} }
and $\mu_n$ in~\eqref{md_1} can be replaced by {\em $\textbf{E}X_1$}, without affecting the result.
\end{theorem}

We now turn to the statement of our results on moderate deviations for $\tn$ in the case of heavy trimming.
Define four numbers
\begin{equation*}
\label{c_an_h}
\begin{split}
a_1&= \liminf_{\nty}\, \an , \qquad  \qquad a_2= \limsup_{\nty}\, \an , \\
b_1&= \liminf_{\nty}\, (1-\bn), \qquad  b_2= \limsup_{\nty}\,
(1-\bn) .
\end{split}
\end{equation*}
Now, we will assume that
\begin{equation}
\label{abc}
0< a_1, \quad b_2 < 1 \quad \text{and}\quad   a_2<b_1.
\end{equation}
In this case no moment assumptions are  needed  for
the asymptotic normality of $\tn$ and for related properties, whereas some
smoothness of  $F^{-1}$ at the points where trimming occurs becomes essential (see, e.g., \cite{cso_hae_mas:1988, gri_helm:2014}, see also the discussion after  the statement of  Theorem~\ref{thm1}).

Let us introduce two sequences of the auxiliary functions:
\begin{equation}
\label{g_h}
\begin{split}
G_n(t)= &F^{-1}\Bigl(\an+t\sqrt{\frac{\an\log n}{n}}\Bigr)-F^{-1}\bigl(\an\bigr),\\
H_n(t)= &F^{-1}\Bigl(1-\bn+t\sqrt{\frac{\bn\log n}{n}}\Bigr)-F^{-1}\bigl(1-\bn\bigr),
\end{split}
\end{equation}
$t\in\mathbb{R}$. Note that, for a~fixed $t$, we have $\an+t\sqrt{\frac{\an\log n}{n}}=\an\bigl(1+t\sqrt{\frac{\log n}{\kn}}\bigr)=\an (1+o(1))$ as $\nty$ (by~\eqref{abc}) and $1-\bn+t\sqrt{\frac{\bn\log n}{n}}=1-\bn \bigl(1+o(1)\bigr)$. In particular this implies that the functions introduced in~\eqref{g_h} are well-defined for all sufficiently large~$n$.

Now we are in a~position to state our first result for the heavy trimmed means.
\begin{theorem}
\label{thm3}
Suppose that condition~\eqref{abc} is satisfied.  In addition, assume
that there exists $\eps>0$ such that  for each \ $t\in \mathbb{R}$
\begin{equation}
\label{cgh}
G_n(t)=O\bigl((\log n)^{-(1+\eps)}\bigr),\quad H_n(t)=O\bigl((\log n)^{-(1+\eps)}\bigr),
\end{equation}
as $\nty$. Then, for each  $c>0$ and  $A>0$,  relations~\eqref{md_1} hold true uniformly in the range
$-A\leq x\leq c\sqrt{\log n}$.
\end{theorem}

\begin{remark} {\em
We notice that $G_n$ differs from the difference in the first line in~\eqref{psi} only by the presence of the logarithm under the root sign, the same remark can be applied to~$H_n$. Thus,
one can see that condition~\eqref{cgh} is somewhat stronger than the $iff$ condition of asymptotic normality by  \cite{cso_hae_mas:1988} (cf.~\eqref{q_0}), but it  enables us to obtain the results on moderate deviations for~$\tn$.
}
\end{remark}
\begin{remark} {\em
It is also worth noting that~\eqref{cgh} holds true if, for instance, the inversion $F^{-1}$ satisfies a~H\"{o}lder condition of degree~$d$ (for some $d>0$) in an~open set containing all limit points of the sequences $\an$ and $1-\bn$. However, the H\"{o}lder condition would be excessive for our present purposes, as  it  can provide us with the stronger results on large deviations (i.e. the deviations in the range of the form $-A\leq x \leq o(n^{r_d})$, for some $0\leq r_d\leq 1/6$) in the case of heavy trimming (cf.~\cite{gri:2016}). 
}
\end{remark}

Finally, we state a~result on  moderate deviations,  parallel to Theorem~\ref{thm2}, but  now for the case of heavy trimming. 
A~very mild moment condition  will be required now to ensure the existence of the variance of $\tn$.

\begin{theorem}
\label{thm4}
Suppose that the conditions of Theorem~\ref{thm3} hold true. 
In addition, assume that {\em $\textbf{E}|X_1|^{\gamma}<\infty$ } for some $\gamma>0$.
Then {\em
\begin{equation}
\label{t4_1}
n^{1/2} (\textbf{E}\tn-\mu_n)=O\Bigl((\log n)^{-(1+\eps)}\Bigr), \qquad\qquad
\end{equation}
\begin{equation}
\label{t4_3}
\frac{\sqrt{ {Var}(\tn)}}{\siw/\sqrt{n}}= 1+ O\Bigl((\log n)^{-(1+\eps)}\Bigr),\
\end{equation}
}
as \,$\nty$, where $\eps$ is as in~\eqref{cgh}.

Moreover,  relations~\eqref{md_1} remain valid, uniformly in the range
$-A\leq x\leq c\sqrt{\log n}$, for each $c>0$ and $A>0$, if we replace $\mu_n$ and/or $\siw$ in it  by {\em $\textbf{E}\tn$ } and {\em $\sqrt{n {Var}(\tn)}$} respectively.
\end{theorem}

\section{On our approach}
\label{lemmas}

Define a binomial r.v. $N_\nu= \sharp \{i : X_i \le \xi_{\nu}\}$, $0<\nu<1$. Set $A_n=F_n(\qa)=N_{\an}/n$, $B_n=1-F_n(\qb)=(n-N_{1-\bn})/n$.  Consider the mean $\overline{W}_n=\frac 1n\sum\limits_{i=1}^n\wi$ of i.i.d. Winsorized r.v.'s  defined in~\eqref{wi}.

The next lemma will be  crucial for our proofs as it provides us with
an~approximation of $T_n$ by sums of i.i.d. r.v.'s
\begin{lemma}
\label{lem1}
The following representation is valid {\em
\begin{equation}
\label{rep}
\tn-\mu_n=\overline{W}_n-\textbf{E}\overline{W}_n+R_n,
\end{equation}
} where
\begin{equation}
\label{rep1}
R_n=\int_{\an}^{A_n}[F_n^{-1}(u)-F^{-1}(\an)]\, du-\int_{1-\bn}^{1-B_n}[F_n^{-1}(u)-F^{-1}(1-\bn)]\, du.
\end{equation}
\end{lemma}

In fact,  representation~\eqref{rep} follows from the first part of the proof of Lemma~2.1 in     \cite[cf.~relation~(2.13)]{gri_helm:2014}. 
For the convenience of the reader we present  briefly its proof. 
%

\medskip
\noindent{\bf Proof.}  Let $\win^{(n)}$ be the order statistics corresponding to the sample $\wi$, $i=1,\dots,n$. Since
$\win^{(n)}=\xin$ for $N_{\an}+1\leq i\leq N_{1-\bn}$, we can write
\begin{equation*}
\label{l11}
\overline{W}_n = \frac 1n\sum\limits_{i=1}^n\wi=A_n\qa+\frac 1n\sum\limits_{i=N_{\an}+1}^{N_{1-\bn}}\xin +
B_n\qb.
\end{equation*}
Using the latter formula, as a result of the direct and simple computations, we obtain
\begin{equation}
\label{l12}
\begin{split}
&T_n-\mu_n -[\overline{W}_n-\textbf{E}\overline{W}_n]=\frac 1n\Bigl[ sgn(N_{\an}-\, \kn)\sum_{i=(\kn\wedge
N_{\an})+1}^{N_{\an}\vee
\kn}(X_{i:n}-\qa)\\
-&sgn(N_{1-\bn}-(n-\mn))\sum_{i=((n-\mn)\wedge
N_{1-\bn})+1}^{N_{1-\bn}\vee \, (n-\mn)}(X_{i:n}-\qb)\Bigr],
\end{split}
\end{equation}
where $sgn(s)=s/|s|$, $sgn(0)=0$. The r.h.s. of~\eqref{l12} is equal to $R_n$. The lemma is proved.
\begin{remark} {\em
It is worth  noting  that  under additional smoothness conditions on $F^{-1}$ (in an~open set containing  the limit points of the sequences $\an$ and $1-\bn$)  representation~\eqref{rep} can be extended to a~$U$-statistic type approximation for $\tn$. We refer to Lemma~2.1 in \cite{gri_helm:2014} for the details. In order to get the quadratic term of the $U$-statistic approximation, in the cited paper we apply some special Bahadur--Kiefer representations of von Mises  statistic type for intermediate sample quantiles obtained in \cite{gri_helm:2012}.
It should be also noted that the idea to approximate the trimmed sums by the sums of i.i.d. Winsorized r.v.'s (as a~linear term of the approximation) plus a~quadratic $U$-statistic term based on the Bahadur type representations was first proposed and implemented  in Gribkova and Helmers~(2006,~2007), where the validity of the one-term Edgeworth
expansions  for a (Studentized) trimmed mean and its bootstrapped version was proved and simple explicit formulas of  these expansions were found.
This approach can be  extended to the case of  trimmed $L$-statistics (cf. \cite{gri:2016}).
 }
\end{remark}

The following lemma on bounds for absolute moments of order statistics will be applied in the proof of Theorems~\ref{thm1}-\ref{thm4}.  
This lemma was obtained in
 \cite[Theorem~1]{gri:1995}, so  we present  here only its statement.

Let $k$ and $\delta$ be arbitrary positive numbers . Put $\rho=k/\delta$ and set $\al_i=i/(n+1)$.
\begin{lemma}
 \label{lem2} {\em(\cite{gri:1995})} For all $n\geq 2\rho+1$ and for all $i$ such that $\rho \leq i\leq n-\rho+1$ the following inequality holds {\em
\begin{equation}
\label{bound}
\textbf{E}|\xin|^k <C(\rho)\lf [\al_i(1-\al_i)]^{-1}\textbf{E}|X_1|^{\delta} \rf^{\rho},
\end{equation} }
where one can put $C(\rho)=2\sqrt{\rho}\exp(\rho+7/6)$.
\end{lemma}
Obviously, estimate~\eqref{bound} is  non-trivial only when  $\textbf{E}|X_1|^{\delta}<\infty$. We
also emphasize the fact that the case $k>\delta$, the most interesting  for us here, is allowed in~\eqref{bound}.
\section{Proofs}
\label{proof}
In this section we prove Theorems~\ref{thm1}-\ref{thm4} stated in Section~\ref{imtro}.

\medskip
\noindent{\bf Proof of Theorem~\ref{thm1}}. 
It suffices to prove the first relation in~\eqref{md_1} (the second relation follows from
the first one if we replace $X_i$ by $-X_i$). Define the $df$'s
\begin{equation}
\label{d.f.s}
F_{\tn}(x)=\textbf{P}\Bigl(\frac{\sqrt{n}(\tn-\mu_n)}{\siw}<x\Bigr);\quad F_{\wn}(x)=\textbf{P}\Bigl(\frac{\sqrt{n}(\overline{W}_n-\textbf{E}\overline{W}_n)}{\siw}<x\Bigr).
\end{equation}
 Applying the classical Slutsky argument to $\tn-\mu_n=\overline{W}_n-\textbf{E}\overline{W}_n+R_n$ (cf.~\eqref{rep})
 gives, for $\dn >0$, that $1-F_{\tn}(x)$ is bounded from above and below by respectively
 \begin{equation*}
\label{p1}
1-F_{\wn}(x-\dn) +\textbf{P}(n^{1/2}\siw^{-1}|R_n|>\dn)
\end{equation*}
and
\begin{equation*}
\label{p2}
1-F_{\wn}(x+\dn) -\textbf{P}(n^{1/2}\siw^{-1}|R_n|>\dn).
\end{equation*}
Now we choose $\dn=[\log(1+n)]^{-d}$ with $d=\gamma\frac{p-2}{2p}-\frac 12$, where $\gamma$ is the constant from condition~\eqref{c_an}. Observe that $d>1/2$ (cf.~\eqref{c_an}) and hence $\dn\sqrt{\log n}\to 0$ as $\nty$.

Thus in order to prove our theorem, it suffices to show that
\begin{equation}
\label{p3}
1-F_{\wn}(x \pm \dn)=[1-\Phi(x)](1+o(1))
\end{equation}
and
\begin{equation}
\label{p4}
\textbf{P}(n^{1/2}\siw^{-1}|R_n|>\dn)=[1-\Phi(x)]o(1),
\end{equation}
as $\nty$, uniformly in the range  $-A\leq x\leq c\sqrt{\log n}$ \  $(A>0)$. 
Fix an~arbitrary $A>0$.

Let us  prove~\eqref{p3}. Define  $Z_i^{(n)}=\siw^{-1}(\wi-\textbf{E}\wi)$, $i=1,\dots,n$, $n=1,2,\dots$, and note  that these r.v.'s form a~triangular series of i.i.d. in each series r.v.'s. Further,  $\textbf{E}Z_i^{(n)}=0$, ${Var}Z_i^{(n)}=1$, and  $\textbf{E}|Z_i^{(n)}|^p $ is bounded from above uniformly in $n$ (because $\wi$ is
the Winsorized $X_i$, and $\textbf{E}|X_1|^p<\infty $). An~application of a~now classical result by \cite{rubin_set:1965} for scheme series yields that
\begin{equation}
\label{p5}
1-F_{\wn}(x \pm \dn)=[1-\Phi(x\pm \dn)](1+o(1))
\end{equation}
uniformly for $x$ in the range $0\leq x \pm \dn\leq c_1 \sqrt{\log n}$, for every $c_1$ such that $p>c_1^2+2$, where  we choose $c_1>c$. Furthermore, by CLT for $Z_i^{(n)}$, $i=1,\dots,n$, $n=1,2,\dots$, relation~\eqref{p5} holds true  uniformly in the range $-A_1\leq x\pm \dn \leq 0$ for each $A_1>0$. Set $q=\sup_{n\in\mathbb{N}}|\dn|=(\log 2)^{-d}$ and put $A_1=A+q$. Then we obtain the validity of~\eqref{p5} uniformly in the range $-A\leq x\leq c \sqrt{\log n}+r_n$,  where $r_n=(c_1-c)\sqrt{\log n}-q\to +\infty$ as $\nty$. At this point we apply Lemma~A1 of \cite{vand_verav:1982},   where  the asymptotic property of~$\Phi$ we need here is given in a~convenient form.
Since $\dn \sqrt{\log n}=o(1)$, due to that lemma we obtain that  $1-\Phi(x\pm \dn)=[1-\Phi(x)](1+o(1))$, uniformly in the range $-A\leq x\leq c \sqrt{\log n}$.
Combining~\eqref{p5} and the latter arguments, we find that~\eqref{p3} is valid  uniformly in the required range.

Let us prove~\eqref{p4}. We first write
\begin{equation}
\label{p6}
\textbf{P}(n^{1/2}\siw^{-1}|R_n|>\dn)\leq \textbf{P}(n^{1/2}\siw^{-1}|R_{n,\al}|>\dn/2)
+ \textbf{P}(n^{1/2}\siw^{-1}|R_{n,\be}|>\dn/2),
\end{equation}
where $R_{n,\al}$, $R_{n,\be}$ denote  the first and second integrals in~\eqref{rep1} respectively. We will  prove~\eqref{p4} for the first probability  on the r.h.s. of~\eqref{p6}, the treatment for the second one is similar and therefore omitted. In view of our moment assumption, $\siw\to\si=[{Var}(X_1)]^{1/2}$, so  it is sufficient  to prove that
\begin{equation}
\label{p7}
\textbf{P}(n^{1/2}|R_{n,\al}|>L\dn)=[1-\Phi(x)]o(1),
\end{equation}
as $\nty$, uniformly in the range  $-A\leq x\leq c\sqrt{\log n}$, where $L$ stands for  a~positive constant not  depending on $n$, which may change  its value  from line to line.
Notice that
\begin{equation}
\label{p8}
\frac 1{1-\Phi(x)}\leq \frac 1{1-\Phi(c\sqrt{\log n})} \sim c\sqrt{\log n}\, n^{c^2/2},
\end{equation}
for $x\in [-A,c\sqrt{\log n}]$. Hence~\eqref{p7} is implied by
\begin{equation}
\label{p9}
\textbf{P}(n^{1/2}|R_{n,\al}|>L\dn)=o\lr(\log n)^{-1/2} n^{-c^2/2} \rr.
\end{equation}
Let us prove~\eqref{p9}. We have $n^{1/2}|R_{n,\al}|=\bigl| n^{-1/2}\sum_{i=(\kn\wedge
N_{\an})+1}^{N_{\an}\vee
\kn}(\xin-\qa)\bigr|$, and by the monotonicity in $i$ of the difference $\xin-\qa$ we obtain
\begin{equation}
\label{p10}
n^{1/2}|R_{n,\al}| \leq n^{-1/2} |N_{\an}-\kn|\,|\xkn-\qa|.
\end{equation}
Let $U_1,\dots,U_n$ be a~sample of independent $(0,1)$-uniform distributed r.v.'s, and let $\uin$ denote the corresponding order statistics. Set $\mau=\sharp \{i : U_i \le \an \}$.  Since the joint
distribution of  $\xin$ and $N_{\an}$ coincides with the joint distribution of $F^{-1}(\uin)$ and $\mau$, $i=1,\dots,n$, we have
\begin{equation}
\label{p11}
\textbf{P}(n^{1/2}|R_{n,\al}|>L\dn)\leq \textbf{P}(n^{-1/2} |\mau-\kn||F^{-1}(\ukn)-F^{-1}(\an)|>L\dn)\leq P_1 +P_2,
\end{equation}
where
\begin{equation*}
\label{p12}
\begin{split}
P_1=&\textbf{P}\bigl(|\mau-\kn|>c_1\sqrt{\kn\log n}\bigr),\\
P_2=&\textbf{P}\bigl(n^{-1/2} \sqrt{\kn\log n}|F^{-1}(\ukn)-F^{-1}(\an)|>L\dn\bigr),
\end{split}
\end{equation*}
and  $c_1$ is as before, i.e., $2+c_1^1<p$, \ $c<c_1$. For $P_1$ by Bernstein's inequality we obtain
\begin{equation}
\label{p13}
P_1 \leq 2 exp(-h_n),
\end{equation}
with
\begin{equation*}
\label{p14}
h_n=
\frac{c_1^2\log n}{2[1-\an+(c_1/3)\sqrt{\log n/\kn}\,\,(\an\vee 1-\an)]}.
\end{equation*}
Since $\an\to 0$ and $\log n/\kn\to 0$, $\nty$, we get that  $h_n\sim \frac{c_1^2\log n}{2}$. Hence   $h_n>\frac{c_2^2\log n}{2}$, for some $c_2$ such that $c<c_2<c_1$ and
all sufficiently large $n$. Then, relations~\eqref{p8} and~\eqref{p13} together yield
\begin{equation}
\label{p15}
P_1 \leq 2 n^{-c_2^2/2}=[1-\Phi(x)]o(1),
\end{equation}	
uniformly in the range $-A\leq x\leq c\sqrt{\log n}$.

It remains to estimate $P_2$ on the r.h.s. in~\eqref{p11}. We have
\begin{equation}
\label{p16}
\begin{split}
P_2=&\textbf{P}\bigl(\an^{1/2} |F^{-1}(\ukn)-F^{-1}(\an)|>L\dn/\sqrt{\log n}\bigr) \\
\leq & \textbf{P}\bigl(\an^{1/2} |F^{-1}(\ukn)|>L\dn/\sqrt{\log n}-\an^{1/2}|F^{-1}(\an)|\bigr) \\
=&\textbf{P}\bigl(\an^{1/2} |F^{-1}(\ukn)|>L\dn/\sqrt{\log n}\bigl[1-\sqrt{\log n}\,\,\dn^{-1}\,\an^{1/2}|F^{-1}(\an)|\bigr]\bigr).
\end{split}
\end{equation}
Since $\textbf{E}|F^{-1}(U_1)|^p<\infty$, we have $\an |F^{-1}(\an)|^p\to 0$, $\nty$, and $\an^{1/2}|F^{-1}(\an)|
=o\bigl(\an^\frac{p-2}{2p}\bigr)$. From the other hand,  $\an=O\bigl((\log n)^{-\gamma}\bigr)$ (due to condition~\eqref{c_an}), and  $\sqrt{\log n}\,\,\dn^{-1}\leq (\log (1+n))^{\gamma\frac{p-2}{2p}}$  (by the choice of $\dn$). The latter computations yield that the second term within square brackets on the r.h.s. in~\eqref{p16} is $o(1)$, and hence it can be omitted. Set $p_n=\textbf{E}\ukn=\frac{\kn}{n+1}$, define $\mathbb{V}_n(p_n)=\sqrt{n}(\ukn-p_n)$, and let $c_1$ be as before (i.e., $2+c_1<p$, \ $c<c_1$).
Then we should evaluate
\begin{equation}
\label{p17}
\textbf{P}\bigl(\an^{1/2} |F^{-1}(\ukn)|>L\dn/\sqrt{\log n}\bigr)\leq P_3 + P_4,
\end{equation}
\begin{equation*}
\label{p18}
\begin{split}
 P_3 =&\textbf{P}\bigl(\big\{\an^{1/2} |F^{-1}(\ukn)|>L\dn/\sqrt{\log n}\big\}\bigcap \big\{|\mathbb{V}_n(p_n)|\leq c_1\sqrt{\an\log n} \big\}\bigr),\\
 P_4 =&\textbf{P}\bigl( |\mathbb{V}_n(p_n)|\geq c_1\sqrt{\an\log n} \bigr).
\end{split}
\end{equation*}

In order to estimate $P_4$, we can apply Inequality~1 from  \cite[page~453]{shor_weln:1986}. Then we obtain
\begin{equation}
\label{p19}
P_4\leq \exp\Bigl[ -c_1^2\frac{\an\log n}{2p_n}\Bigr]+\exp\Bigl[ -c_1^2\frac{\an\log n}{2p_n}\widetilde\psi(t_n)\Bigr],
\end{equation}
where $\widetilde\psi$ is the function defined in~\cite[page~453, formula~(2)]{shor_weln:1986},
$t_n=c_1\frac{\sqrt{\an\log n}}{p_n\sqrt{n}}=c_1\sqrt{\frac{n+1}{n}}\sqrt{\frac{\log n}{\kn}}$. By condition~\eqref{c_kn}, $t_n\to 0$ as $\nty$, hence $t_n>-1$ for all sufficiently large $n$, and by Proposition~1 in~\cite[page~455, relation~(12)]{shor_weln:1986}, we find that $\widetilde\psi(t_n)\geq \frac 1{1+2t_n/3}$. This and relation~\eqref{p19} together yield
\begin{equation}
\label{p20}
P_4\leq 2\exp\Bigl[ -c_2^2\frac{\log n}{2}\Bigr] =2n^{-c_2^2/2},
\end{equation}
for each $c_2$ such that $c<c_2<c_1$ and for all sufficiently large $n$. Let $M_n$ denote $F^{-1}\bigl( p_n-c_1\frac{\sqrt{\an\log n}}{\sqrt{n}}\bigr)\vee F^{-1}\bigl( p_n+c_1\frac{\sqrt{\an\log n}}{\sqrt{n}}\bigr)$, then by monotonicity of $F^{-1}$, we get
\begin{equation}
\label{p21}
P_3\leq \textbf{P}\bigl(\an^{1/2}M_n>L\dn/\sqrt{\log n} \bigr).
\end{equation}
Observe that $p_n \pm c_1\frac{\sqrt{\an\log n}}{\sqrt{n}}=\an\bigl(1 \pm O \bigl(\sqrt{\frac{\log n}{\kn}} \bigr)\bigr)=\an(1+o(1))$, and using our moment assumption similarly as before, we find that $\an^{1/2}M_n=o\bigl( (\log n)^{-\gamma\frac{p-2}{2p}}\bigr)$, whereas $\dn/\sqrt{\log n}\geq (\log (1+n))^{-\gamma\frac{p-2}{2p}}$. Hence  the quantity on the r.h.s. in~\eqref{p21} is equal to zero for all sufficiently large $n$. Relations~\eqref{p11},~\eqref{p15}-\eqref{p17} and \eqref{p20}-\eqref{p21} together imply~\eqref{p9}, which  entails~\eqref{p4}. The theorem is proved.

\medskip
\noindent{\bf Proof of Theorem~\ref{thm2}}. Let us first prove~\eqref{t2_1}. By Lemma~\ref{lem1}, we get  $\textbf{E}\tn-\mu_n=\textbf{E}R_{n,\al} -\textbf{E}R_{n,\be}$. Here and later on,  we keep the notation introduced in the proof of the previous theorem. 
Then,
\begin{equation}
\label{p22}
n^{1/2}|\textbf{E}\tn-\mu_n|\leq n^{1/2}\textbf{E}|R_{n,\al}| +n^{1/2}\textbf{E}|R_{n,\be}|.
\end{equation}
We will estimate only the first term on the r.h.s. in~\eqref{p22} (obviously the handling for the second one is  similar). As in the proof of Theorem~\ref{thm1} (cf.~\eqref{p10}-\eqref{p11}), we find that
\begin{equation}
\label{p23}
\begin{split}
n^{1/2}\textbf{E}|R_{n,\al}| \leq &n^{-1/2}\textbf{E}(|\mau-\kn||F^{-1}(\ukn)-F^{-1}(\an)|)\\
\leq & n^{-1/2}\bigl[ \textbf{E}(\mau-\kn)^2 \textbf{E}(F^{-1}(\ukn)-F^{-1}(\an))^2 \bigr]^{1/2}\\
=& \an^{1/2}(1-\an)^{1/2} \bigl[ \textbf{E}(F^{-1}(\ukn)-F^{-1}(\an))^2 \bigr]^{1/2}\\
\leq & \an^{1/2}(1-\an)^{1/2} 2^{1/2}\bigl[( \textbf{E}(F^{-1}(\ukn))^2 )^{1/2}+ |F^{-1}(\an))|  \bigr].
\end{split}
\end{equation}
By our moment assumption, we have $|F^{-1}(\an))|=o(\an^{-1/p})$. In order to estimate $( \textbf{E}(F^{-1}(\ukn))^2 )^{1/2}$ on the r.h.s. in~\eqref{p23}, we apply Lemma~\ref{lem2}. Then we obtain
\begin{equation*}
\label{p24}
(F^{-1}(\ukn))^2 )^{1/2}\leq (\an (1-\an))^{-1/p} (\textbf{E}|X_1|^p)^{1/p}=O(\an^{-1/p}).
\end{equation*}
Hence, the quantity on the r.h.s. in~\eqref{p23} is of the order $O(\an^{\frac{p-2}{2p}})=O\bigl((\log n)^{-\gamma \frac{p-2}{2p}}\bigr)=o\bigl((\log n)^{-1}\bigr)$. Thus~\eqref{t2_1} follows.

Next we prove~\eqref{t2_4} (using the additional condition~\eqref{c_an2}). We write
\begin{equation}
\label{p25}
\begin{split}
n^{1/2}|\textbf{E}X_1-\mu_n|=& n^{1/2}\Bigl|\int_0^{\an}F^{-1}(u)\,du + \int_{1-\bn}^1 F^{-1}(u)\,du\Bigr|\\
\leq &n^{1/2}\Bigl(\int_0^{\an}|F^{-1}(u)|\,du \Bigr)+ n^{1/2}\Bigl(\int_{1-\bn}^1 |F^{-1}(u)|\,du\Bigr).
\end{split}
\end{equation}
As before, we estimate only the first term on the r.h.s. in~\eqref{p25}. Since $\an\to 0$, our moment assumption implies that for every $K>0$ and sufficiently large~$n$
\begin{equation}
\label{p26}
\begin{split}
n^{1/2}\int_0^{\an}|F^{-1}(u)|\,du\leq &Kn^{1/2}\int_0^{\an} u^{-1/p}\,du = Kn^{1/2}\frac{p}{1-p}\,\,\an^{\frac{p-1}{p}}\\
=&o\Bigl(n^{1/2} (n\log n)^{-  1/2}\Bigr)=o\Bigl( (\log n)^{-  1/2}\Bigr)
\end{split}
\end{equation}
(here condition~\eqref{c_an2} was applied). Relations~\eqref{p25} and~\eqref{p26} yields~\eqref{t2_4}.

Let us prove~\eqref{t2_2}. Since under our moment assumption $\siw\to\si$ as $\nty$, it suffices to show that $\si^2-\siw^2=o((\log n)^{-2})$. We have
\begin{equation}
\label{p27}
\si^2-\siw^2=\int\limits_{[0,1]^2}\int\limits_{\setminus [\an, 1-\bn]^2} (u\wedge v - uv)\,dF^{-1}(u)\,dF^{-1}(v).
\end{equation}
Due to the fact that $\an<1/2<1-\bn$ for all sufficiently large $n$ and because of the symmetry in our conditions for  $\an$ and $\bn$, it is sufficient to estimate the integral over the region $0 \leq u\leq \an$, $u\leq v\leq 1/2$. We have
\begin{equation}
\label{p28}
\begin{split}
&\int_0^{\an}\int_u^{1/2} (u\wedge v - uv)\,dF^{-1}(v)\,dF^{-1}(u)
=\int_0^{\an}u\int_u^{1/2} (1 - v)\,dF^{-1}(v)\,dF^{-1}(u)\\
\leq &\int_0^{\an}u\int_u^{1/2}\,dF^{-1}(v)\,dF^{-1}(u)=\int_0^{\an}u [F^{-1}(1/2)- F^{-1}(u)]\,dF^{-1}(u)\\
\leq &|F^{-1}(1/2)| \int_0^{\an}u\,dF^{-1}(u)+ \int_0^{\an}u |F^{-1}(u)|\,dF^{-1}(u).
\end{split}
\end{equation}
For the first integral on the r.h.s. in~\eqref{p28} for sufficiently large $n$ we obtain
\begin{equation*}
\label{p29}
\begin{split}
 \int_0^{\an}u\,dF^{-1}(u)\leq \an\,|F^{-1}(\an)| +K\int_0^{\an}u^{-1/p}\,du\leq K\,\an^{1-1/p}\frac {2p}{p-1},
\end{split}
\end{equation*}
where $K$ is as before. For the second integral on the r.h.s. in~\eqref{p28} similarly we find
\begin{equation*}
\label{p30}
\begin{split}
 \int_0^{\an}u |F^{-1}(u)|\,dF^{-1}(u)\leq K\int_0^{\an}u^{\frac{p-1}{p}}\,dF^{-1}(u)\leq
 K^2\,\an^{1-2/p}\,\,\frac {2p-3}{p-2}.
\end{split}
\end{equation*}

The latter computations imply that the quantity on the r.h.s. in~\eqref{p28} is of the order $O\bigl( \an^{\frac{p-2}{p}}\bigr)$, and by condition~\eqref{c_an} it is $o\bigl( (\log n)^{-2}\bigr)$.

Finally, we prove~\eqref{t2_3}. In fact, we should show that
\begin{equation}
\label{p31_0}
n{Var}(\tn)-\siw^2=o\bigl( (\log n)^{-1}\bigr).
\end{equation}
By Lemma~\ref{lem1},
\begin{equation*}
\label{p31}
\begin{split}
n{Var}(\tn) =n{Var} (\overline{W}_n+R_n)=&\siw^2+n{Var} (R_n)
+2n\, \text{cov}(\overline{W}_n,R_n)\\
\leq &\siw^2+n{Var} (R_n)
+2\siw(n {Var} (R_n) )^{1/2}.
\end{split}
\end{equation*}  
Hence,
\begin{equation}
\label{p32}
n{Var}(\tn)-\siw^2\leq D_n+2\siw D_n^{1/2},
\end{equation}
where
\begin{equation}
\label{p33}
\begin{split}
D_n=n{Var} (R_n)=&n{Var} (R_{n,\al}-R_{n,\be})\\
\leq &n[{Var}(R_{n,\al})+{Var}(R_{n,\be}) +2({Var}(R_{n,\al}){Var}(R_{n,\be}))^{1/2} ].
\end{split}
\end{equation}
It remains to estimate ${Var}(R_{n,\al})$. Write 
\begin{equation}
\label{p34}
\begin{split}
{Var} (R_{n,\al}) \leq \textbf{E}R_{n,\al}^2\leq &n^{-2}\textbf{E}\bigl[(\mau-\kn)^2(F^{-1}(\ukn)-F^{-1}(\an))^2 \bigr]\\
\leq &n^{-2} \Bigl[\textbf{E}(\mau-\kn)^4  \textbf{E}(F^{-1}(\ukn)-F^{-1}(\an))^4 \Bigr]^{1/2}.
\end{split}
\end{equation}
By well-known formula for the forth moment of the binomial r.v., we have $\textbf{E}(\mau-\kn)^4 =n\an(1-\an)[1+3(n-2)\an(1-\an)]< 3n^2\an^2[1/\kn+1]\sim 3n^2\an^2$. Further, we find that
\begin{equation*}
\label{p35}
\textbf{E}(F^{-1}(\ukn)-F^{-1}(\an))^4\leq 8[ \textbf{E}(F^{-1}(\ukn))^4 +(F^{-1}(\an))^4],
\end{equation*}
where by our moment assumption we have $(F^{-1}(\an))^4=o( \an^{-4/p})$ and by Lemma~\ref{lem2},
$\textbf{E}(F^{-1}(\ukn))^4\leq C[\an(1-\an]^{-4/p}\textbf{E}(|X_1|^p)^{4/p}=O( \an^{-4/p})$. Combining~\eqref{p34} and the latter computations, we obtain that ${Var} (R_{n,\al})=O(n^{-1}\an^{1-2/p})$. Similarly, we find that ${Var} (R_{n,\be})O(n^{-1}\bn^{1-2/p})$. Hence, by~\eqref{p33},
\begin{equation}
\label{p36}
D_n=O\Bigl( (\an\vee\bn)^{1-2/p}\Bigr)=o\Bigl( (\log n)^{-2}\Bigr).
\end{equation}
Relations~\eqref{p32} and~\eqref{p36} imply~\eqref{p31_0}. Thus~\eqref{t2_3} follows.

To complete the proof it remains to argue why one can replace the normalizing sequences $\mu_n$, $\siw$ in~\eqref{md_1} as it was stated in  Theorem~\ref{thm2}. In fact, this is implied by Theorem~\ref{thm1}, Lemma~A.1 of \cite{vand_verav:1982} and relations~\eqref{t2_1}-\eqref{t2_3},~\eqref{t2_4}. Indeed, let $A_n$ denote $\textbf{E}\tn$ or -- under additional assumption~\eqref{c_an2} --  $\textbf{E}X_1$, and let $B_n$ denote  $\sqrt{n{Var}(\tn)}$ or   $\si$. Put $\lambda_n=B_n/\siw$, $\nu_n=n^{1/2}\siw^{-1}(A_n-\mu_n)$. By~\eqref{t2_1}-\eqref{t2_3} and~\eqref{t2_4}, in each of these options we have $\nu_n=o\bigl((\log n)^{-1/2}\bigr)$ and $\lambda_n-1=o\bigl((\log n)^{-1}\bigr)$. Fix an~arbitrary $A>0$ and $c_1$ such that $c<c_1<\sqrt{p-2}$.
Take $A_1>A$.  By Theorem~\ref{thm1},
\begin{equation}
\label{p37}
\begin{split}
\textbf{P}\bigl(B_n^{-1}n^{1/2}(\tn-A_n) >x \bigr)=&\textbf{P}\bigl(\siw^{-1}n^{1/2}(\tn-\mu_n)>\lambda_n x+\nu_n \bigr)\\
=&[1-\Phi(\lambda_n x+\nu_n)](1+o(1)),
\end{split}
\end{equation}
for $x$ such that $-A_1\leq \lambda_n x+\nu_n \leq c_1\sqrt{\log n}$. Note that for all sufficiently large $n$ we have   $\lambda_n^{-1}(A_1-\nu_n)>A$ and $\lambda_n^{-1}(c_1\sqrt{\log n}-\nu_n)> c\sqrt{\log n}$. Hence~\eqref{p37} holds  uniformly in the range
$-A \leq x\leq c\sqrt{\log n}$. Finally, we note that since $\sqrt{\log n}\Bigl(|\lambda_n-1|^{1/2}\vee |\nu_n|\Bigr)\to 0$ as $\nty$, the Lemma~A.1 of \cite{vand_verav:1982}  implies that $[1-\Phi(\lambda_n x+\nu_n)]=[1-\Phi(x)](1+o(1))$, uniformly in the range
$-A \leq x\leq c\sqrt{\log n}$. This fact and relation~\eqref{p37} together yield the result desired.  The theorem is proved.

\medskip
\noindent{\bf Proof of Theorem~\ref{thm3}}. Similarly as in the proof of Theorem~\ref{thm1}, we start with the application of the Slutsky argument, where now we set $\dn=(\log (1+n))^{-(1/2+\eps_1)}$ with an arbitrary $\eps_1$ such that  $0<\eps_1<\eps$, where  $\eps$ is as in~\eqref{cgh}. Then, we notice again  that it suffices to prove that~\eqref{p3} and~\eqref{p4} are valid  uniformly in the range  $-A\leq x\leq c\sqrt{\log n}$ ($A>0$) (but now for each~$c>0$).  Fix arbitrary $c,\,A>0$.

First we prove the validity of~\eqref{p3}. Since $0<a_1$ and $b_1<1$, the r.v.'s $|\wi|$ are bounded from above  uniformly in $n$. Hence the Cram\'{e}r condition for r.v.'s $\wi$ is satisfied  uniformly in $n$, i.e. $\textbf{E}e^{h|\wi|}\leq M<\infty$, for some positive $M$ and all $h>0$ and  $n\in\mathbb{N}$. Then an~application of a~large deviations result for the sum of i.i.d. r.v.'s $\overline{W}_n$ (cf.,e.g.,~\cite{feller:1943, petrov:1975}) yields
\begin{equation}
\label{p40}
1-F_{\wn}(x \pm \dn)=[1-\Phi(x\pm \dn)](1+o(1)),
\end{equation}
as $\nty$, uniformly in the range $-A_1\leq x \pm \dn =o(n^{1/6})$. Put  $A_1=A+\sup_{n\in \mathbb{N}}\dn=A+(\log 2)^{-(1/2+\eps_1)}$ and notice that $c\sqrt{\log n}=o(n^{1/6} -\dn)$. Hence we obtain that~\eqref{p40} is valid uniformly in the range  $-A\leq x\leq c\sqrt{\log n}$.
An~application of Lemma~A.1 by \cite{vand_verav:1982} yields that  $[1-\Phi(x\pm \dn)]=[1-\Phi(x)](1+o(1))$  uniformly in the range  $-A\leq x\leq c\sqrt{\log n}$.

Let us prove that~\eqref{p4} is valid uniformly in the range  $-A\leq x\leq c\sqrt{\log n}$. As before, we  prove it only for $R_{n\al}$, i.e. for the first part  of $R_n$. Since $0<a_1$, \ $b_1<1$ and  because $\liminf_{\nty}\siw>0$, the variance $\siw^2={Var}(\wi)$  is bounded away  from zero and infinity. So, taking into account~\eqref{p8} (cf. also~\eqref{p9}), we see that it suffices to prove that 
\begin{equation}
\label{p41}
\textbf{P}(n^{1/2}|R_{n,\al}|>L\dn)=o\lr(\log n)^{-1/2} n^{-c^2/2} \rr,
\end{equation}
where $L$ denote a~positive constant not depending on $n$, which may change  its value.
Similarly as in the proof of Theorem~\ref{thm1}, we write
\begin{equation}
\label{p42}
\begin{split}
&\textbf{P}(n^{1/2}|R_{n,\al}|>L\dn)\leq \textbf{P}(n^{-1/2}|\mau-\kn||F^{-1}(\ukn)-F^{-1}(\an)|>L\dn)\\
\leq &\textbf{P}(|\mau-\kn|>c_1\sqrt{\an n\log n}) + \textbf{P}(\sqrt{\log(n+1)}|F^{-1}(\ukn)-F^{-1}(\an)|>L\dn),
\end{split}
\end{equation}
where $c_1>c$. For the first probability on the r.h.s. in~\eqref{p42} we find (cf.~\eqref{p13}-\eqref{p15}) that for $c<c_2<c_1$ and all sufficiently large $n$ it does not exceed $2 n^{-c_2^2/2}$ which is $[1-\Phi(x)]o(1)$,  as $\nty$, uniformly in the range  $-A\leq x\leq c\sqrt{\log n}$. It remains to evaluate  the second probability on the r.h.s. in~\eqref{p42}. It is equal to
\begin{equation}
\label{p43}
\textbf{P}\bigl(|F^{-1}(\ukn)-F^{-1}(\an)|>\frac{L}{ (\log (1+n))^{(1+\eps_1)}}\bigr) \leq P_1+P_2,
\end{equation}
where
\begin{equation*}
\label{p44}
\begin{split}
P_1=&\textbf{P}\bigl(\big\{|F^{-1}(\ukn)-F^{-1}(\an)|>\frac{L}{ (\log (1+n))^{(1+\eps_1)}} \big\}\bigcap \big\{|\mathbb{V}_n(\an)|\leq c_1\sqrt{\an\log n} \big\}\bigr)\bigr),\\
P_2=&\textbf{P}\bigl( |\mathbb{V}_n(\an)|\geq c_1\sqrt{\an\log n} \bigr),
\end{split}
\end{equation*}
where $c_1$ is as before (i.e., $c_1>c$) and  $\mathbb{V}_n(\an)=\sqrt{n}(\ukn-\an)$. By the monotonicity of $F^{-1}$,  we get
\begin{equation}
\label{p45}
P_1\leq \textbf{P}\bigl(|G_n(c_1)|\vee|G_n(-c_1)|>{L}{ (\log (1+n))^{-(1+\eps_1)}}\bigr),
\end{equation}
and by condition~\eqref{cgh} and due to the fact that  $\eps_1<\eps$, the probability $P_1$ is zero for all sufficiently large~$n$. Finally, by the same way as in the proof of Theorem~\ref{thm1} (cf.~\eqref{p19}-\eqref{p20}),  an~application of Inequality~1 from \cite{shor_weln:1986} yields that  $P_2\leq 2\exp\Bigl[ -c_2^2\frac{\log n}{2}\Bigr] =2n^{-c_2^2/2}$ for all sufficiently large $n$ and some $c_2$ such that $c<c_2<c_1$. The latter computations imply~\eqref{p41}. Hence~\eqref{p4} holds true as $\nty$,  uniformly in the range  $-A\leq x\leq c\sqrt{\log n}$, for each $c>0$ and $A>0$. The theorem is proved.

\medskip
\noindent{\bf Proof of Theorem~\ref{thm4}}. First we prove~\eqref{t4_1},  starting from  relation~\eqref{p22} and noting  that it suffices to evaluate $n^{1/2}\textbf{E}|R_{n,\al}|$. Similarly as in the proof of Theorem~\ref{thm2} (cf.~\eqref{p23}, we first write
\begin{equation}
\label{p46}
\begin{split}
n^{1/2}\textbf{E}|R_{n,\al}| \leq &n^{-1/2}\bigl[ \textbf{E}(\mau-\kn)^2 \textbf{E}(F^{-1}(\ukn)-F^{-1}(\an))^2 \bigr]^{1/2}\\
=& \an^{1/2}(1-\an)^{1/2} \bigl[ \textbf{E}(F^{-1}(\ukn)-F^{-1}(\an))^2 \bigr]^{1/2}\\
\leq &K \bigl[ \textbf{E}(F^{-1}(\ukn)-F^{-1}(\an))^2 \bigr]^{1/2},
\end{split}
\end{equation}
where $K$ is some positive constant not depending on $n$. Set $\mathcal{E}=\{|\ukn-\an|\leq c\sqrt{\frac{\an\log n}{n}}\}$ and let $\textbf{1}_{\mathcal{E}}$ denote  the indicator of the event~$\mathcal{E}$. Then we can write
\begin{equation}
\label{p47}
\begin{split}
&\bigl[ \textbf{E}(F^{-1}(\ukn)-F^{-1}(\an))^2 \bigr]^{1/2}\\
=&\bigl[ \textbf{E}\bigl((F^{-1}(\ukn)-F^{-1}(\an))^2\textbf{1}_{\mathcal{E}}\bigr)+ \textbf{E}\bigl((F^{-1}(\ukn)-F^{-1}(\an))^2\textbf{1}_{\overline{\mathcal{E}}}\bigr)\bigr]^{1/2}\\
\leq & [(G_n^2(c)\vee G_n^2(-c)]^{1/2} +\bigl[\textbf{E}\bigl((F^{-1}(\ukn)-F^{-1}(\an))^4 \bigr]^{1/4}
  \bigl[ \textbf{P}(\overline{\mathcal{E}}) \bigr]^{1/4}.
\end{split}
\end{equation}
By condition~\eqref{cgh}, the first term on the r.h.s. in~\eqref{p47} is of the order $O\bigl((\log n)^{-(1+\eps)} \bigr)$. By our moment assumption and Lemma~\ref{lem2}, the first factor of the second term  on the r.h.s. in~\eqref{p47} is $O(1)$, and by the Inequality~1 from   \cite{shor_weln:1986}, for the second factor
we get $\textbf{P}(\overline{\mathcal{E}})=o(n^{-c^2/8})$. Hence, the second term on the r.h.s. in~\eqref{p47} is of negligible order for our purposes and contributes to the first one. The latter estimates and~\eqref{p46}-\eqref{p47} imply~\eqref{t4_1}.

Finally, we prove~\eqref{t4_3}. Similarly to the proof of Theorem~\ref{thm2} (cf.~\eqref{p31_0}), we notice that it suffices to show that
\begin{equation}
\label{p48}
n{Var}(\tn)-\siw^2=O\bigl( (\log n)^{-(1+\eps)}\bigr).
\end{equation}
Then, we repeat relations~\eqref{p32}-\eqref{p34} from the proof of Theorem~\ref{thm2}. Thus, we see that one should evaluate  the quantity on the r.h.s. in~\eqref{p34}, but now under the conditions of Theorem~\ref{thm4}. Similarly as before, we find that $\textbf{E}(\mau-\kn)^4 < 3n^2$. So it remains to estimate $[\textbf{E}(F^{-1}(\ukn)-F^{-1}(\an))^4]^{1/2}$,  for which -- by the same way as in~\eqref{p47} -- we get the bound of the order $O\bigl((\log n)^{-2(1+\eps)} \bigr)$.
Hence $D_n=n{Var} (R_n)=O\bigl((\log n)^{-2(1+\eps)} \bigr)$ (cf.~\eqref{p33}), and since $n{Var}(\tn)-\siw^2=O(D_n^{1/2})$, relation~\eqref{p48} follows.

To complete the proof, it remains to argue the possibility of replacing
$\mu_n$ and $\siw$ in relation~\eqref{md_1}  by $\textbf{E}\tn$ and $\sqrt{n{Var}(\tn)}$ respectively. Again we set $A_n=\textbf{E}\tn$, $B_n=\sqrt{n{Var}(\tn)}$,
$\lambda_n=B_n/\siw$, $\nu_n=n^{1/2}\siw^{-1}(A_n-\mu_n)$. Fix arbitrary $c,A>0$, $A_1>A$, $c_1>c$.
Then, using Theorem~\ref{thm3} and the argument below relation~\eqref{p37}, we find that~\eqref{p37} holds  uniformly in the range $-A\leq  x \leq c\sqrt{\log n}$. Furthermore, since  $\sqrt{\log n}\Bigl(|\lambda_n-1|^{1/2}\vee |\nu_n|\Bigr)\to 0$ as $\nty$ (due to~\eqref{t4_1}-\eqref{t4_3}), by Lemma~A.1 of \cite{vand_verav:1982},  we obtain  that $[1-\Phi(\lambda_n x+\nu_n)]=[1-\Phi(x)](1+o(1))$, uniformly in the range $-A \leq x\leq c\sqrt{\log n}$. 
The theorem is proved.

\bibliographystyle{apalike} 
\bibliography{gribkova_ref}

\end{document}